\documentclass[12pt,leqno]{article}
\usepackage{amsfonts}
\usepackage{amsmath, amsthm, amsfonts, amssymb, color}
\usepackage{mathrsfs}
\usepackage{color}
\theoremstyle{definition}

\newcommand{\scr}[1]{\mathscr #1}
\definecolor{wco}{rgb}{0.5,0.2,0.3}

\numberwithin{equation}{section} \theoremstyle{remark}
\newcommand{\ua}{\uparrow}

\title{{\bf  Derivative Formulas in  Measure on Riemannian Manifolds}\footnote{Supported in
 part by  NNSFC (11771326, 11831014, 11921001).} }
\author{
{\bf Panpan Ren$^{b,c)}$,    Feng-Yu Wang$^{a,b)}$  }\\
\footnotesize{$^{a)}$ Center for Applied Mathematics, Tianjin University, Tianjin 300072, China}\\
 \footnotesize{$^{b)}$ Department of Mathematics,
Swansea University, Singleton Park, SA2 8PP, United Kingdom}\\
\footnotesize{$^{c)}$ Mathematical Institute,Woodstock Road, OX2 6GG, University of Oxford}\\
\footnotesize{ 673788@swansea.ac.uk, Panpan.ren@maths.ox.ac.uk; wangfy@tju.edu.cn, F.-Y.Wang@swansea.ac.uk}}
\begin{document}
\allowdisplaybreaks
\def\R{\mathbb R}  \def\ff{\frac} \def\ss{\sqrt} \def\B{\mathbf
B}
\def\N{\mathbb N} \def\kk{\kappa} \def\m{{\bf m}}
\def\ee{\varepsilon}\def\ddd{D^*}
\def\dd{\delta} \def\DD{\Delta} \def\vv{\varepsilon} \def\rr{\rho}
\def\<{\langle} \def\>{\rangle}
  \def\nn{\nabla} \def\pp{\partial} \def\E{\mathbb E}
\def\d{\text{\rm{d}}} \def\bb{\beta} \def\aa{\alpha} \def\D{\scr D}
  \def\si{\sigma} \def\ess{\text{\rm{ess}}}\def\s{{\bf s}}
\def\beg{\begin} \def\beq{\begin{equation}}  \def\F{\scr F}
\def\Ric{\mathcal Ric} \def\Hess{\text{\rm{Hess}}}
\def\e{\text{\rm{e}}} \def\ua{\underline a} \def\OO{\Omega}  \def\oo{\omega}
 \def\tt{\tilde}\def\[{\lfloor} \def\]{\rfloor}
\def\cut{\text{\rm{cut}}} \def\P{\mathbb P} \def\ifn{I_n(f^{\bigotimes n})}
\def\C{\scr C}      \def\aaa{\mathbf{r}}     \def\r{r}
\def\gap{\text{\rm{gap}}} \def\prr{\pi_{{\bf m},\varrho}}  \def\r{\mathbf r}
\def\Z{\mathbb Z} \def\vrr{\varrho} \def\ll{\lambda}
\def\L{\scr L}\def\Tt{\tt} \def\TT{\tt}\def\II{\mathbb I}
\def\i{{\rm in}}\def\Sect{{\rm Sect}}  \def\H{\mathbb H}
\def\M{\mathbb M}\def\Q{\mathbb Q} \def\texto{\text{o}} \def\LL{\Lambda}
\def\Rank{{\rm Rank}} \def\B{\scr B} \def\i{{\rm i}} \def\HR{\hat{\R}^d}
\def\to{\rightarrow}\def\l{\ell}\def\iint{\int}\def\gg{\gamma}
\def\EE{\scr E} \def\W{\mathbb W}
\def\A{\scr A} \def\Lip{{\rm Lip}}\def\S{\mathbb S}
\def\BB{\scr B}\def\Ent{{\rm Ent}} \def\i{{\rm i}}\def\itparallel{{\it\parallel}}
\def\g{{\mathbf g}}\def\Sect{{\mathcal Sec}}\def\T{\mathcal T}\def\BB{{\bf B}}
\def\f{\mathbf f} \def\g{\mathbf g}\def\BL{{\bf L}}  \def\BG{{\mathbb G}}
\def\Bd{{D^E}} \def\BdP{D^E_\phi} \def\Bdd{{\bf \dd}} \def\Bs{{\bf s}} \def\GA{\scr A}
\def\Bg{{\bf g}}  \def\Bdd{{\bf d}} \def\supp{{\rm supp}}\def\div{{\rm div}}
\def\ddiv{{\rm div}}\def\osc{{\bf osc}}\def\1{{\bf 1}}\def\BD{\mathbb D}\def\GG{\Gamma}
\maketitle

\begin{abstract}  We characterise the link of  derivatives in measure, which are  introduced in \cite{AKR,Card,ORS} respectively by different means,
 for functions
on the space $\M$ of finite measures over a Riemannian manifold $M$.    For a reasonable class of functions $f$,
   the extrinsic derivative $D^Ef$ coincides with the linear functional derivative $D^Ff$,
the intrinsic derivative $D^If$ equals to the $L$-derivative $D^Lf$, and
$$D^If(\eta)(x)=  D^{L}f(\eta)(x)= \lim_{s\downarrow 0} \ff 1 s  \nn f(\eta+s \dd_\cdot)(x) = \nn \big\{D^E f (\eta)\big\}(x),     \ \   (x,\eta)\in M\times\M,$$ where $\nn$ is the gradient on $M$,   $\dd_x$ is the Dirac measure
at $x$, and $$D^Ef(\eta)(x):= \lim\limits_{s\downarrow 0} \ff { f(\eta+s \dd_x)-f(\eta)} s,\ \ x\in M$$ is the extrinsic derivative of $f$ at $\eta\in \M$.
This gives a simple way to calculate the intrinsic or $L$-derivative, and is extended to functions of probability measures.
%This provides a simple way to calculate the intrinsic/Lions derivative.
 \end{abstract} \noindent
 \vskip 2cm
\section{Introduction}

To develop analysis on the space of measures, some derivatives in measure have been introduced   by different means,
where the intrinsic and extrinsic derivatives defined  in \cite{AKR,ORS}
  have been used to investigate measure-valued diffusion processes   over Riemannian manifolds (see \cite{KLV,RW19, Sturm,W18,W19} and references therein), and
  the $L$- and linear functional derivatives were investigated in \cite{Card, CDLL} on the Wasserstein space $\scr P_2(\R^d)$ (the the set of all probability measures on $\R^d$ with finite second-order moments).      See  \cite{AMB} and references therein for   calculus and optimal transport on the space of probability measures, and see     \cite{RW19a, Song} for the the Bismut formula and estimates on the $L$-derivative of
distribution dependent SDEs.

In this paper, we aim to  clarify  the link  of   these derivatives, and present formulas for calculations. For a broad range of applications, we will work on the space of finite/probability measures over a Riemannian manifold, which includes $\scr P_2(\R^d)$ as a special example.

 Let $(M,\<\cdot,\cdot\>)$ be a complete  Riemannian manifold,  i.e. $M$ is a differentiable manifold equipped with  the Riemannian metric $\<\cdot,\cdot\>$,
 which is a positive definite smooth bilinear form on the tangent bundle $TM:=\cup_{x\in M}T_xM$ ($T_xM$ is the tangent space at point $x$), such that $M$ is a Polish space under the corresponding Riemannian distance $\rr$.
Let $\M$ denote the class of all nonnegative  finite   measures on $M$ equipped with the weak topology induced by bounded continuous funtions.

For   a fixed point $o\in M$,   let $\rr_o=\rr(o,\cdot)$ be the Riemannian distance function to $o$. Denote $\eta(f)= \int_Mf\d\eta$ for a measure $\eta$ and a function $f\in L^1(\eta)$. For any $p\in [0,\infty)$,  consider   the spaces
$$\M_p:=\big\{\eta\in\M: \eta(\rr_o^p)<\infty\big\},\ \ \scr P_p:= \big\{\eta\in \M_p: \eta(M)=1\big\},\ \ p\in [0,\infty).$$
We will study the above mentioned  derivatives  on $\M_p$ and $\scr P_p.$

For every $p\in [0,\infty)$, $\M_p$ is equipped with the topology that $\eta_n\to\eta$ in $\M_p$ as $n\to\infty$ if and only if
the convergence holds under the weak topology and
$$\lim_{m\to\infty}\sup_{n\ge 1} \eta_n(\rr_o^p 1_{\{\rr_o\ge m\}}) =0.$$
When $p=0,$ this is nothing but the weak topology. When $p>0$, the topology is induced by the $p$-Wasserstein metric
$$\W_p(\gg,\eta):= |\gg(1+\rr_o^p)- \eta(1+\rr_o^p)| + \inf_{\pi\in \C(\gg,\eta)} \big\{\pi(\rr^p) \big\}^{\ff 1 {p\lor 1}},$$
where $\pi\in \C(\gg,\eta)$ means that $\pi$ is a finite measure on $M\times M$ such that
$$\pi(M\times\cdot)= \gg(M)\eta,\ \ \pi(\cdot\times M)= \eta(M) \gg.$$
It is well known that   $(\M_p,\W_p)$ is a Polish space for any $p\in [0,\infty)$.

\

We first recall the extrinsic derivative defined as partial derivative in the direction of Dirac measures, see \cite[Definition 1.2]{ORS}.

\beg{defn}[{\bf Extrinsic derivative}]   Let $p\in [0,\infty)$ and
  $f$ be a real function on $\M_p$.
  \beg{enumerate} \item[$(1)$] $f$  is called extrinsically differentiable on $\M_p$ with derivative $D^Ef$,    if
 $$D^E f(\eta)(x):= \lim_{\vv\downarrow 0} \ff{ f(\eta+\vv\dd_x)-f(\eta)}\vv  \in \R$$  exists  for  all $(x,\eta)\in M\times\M_p.$
 \item[$(2)$]  If    $D^Ef(\eta)(x)$  exists and  is continuous in $(x,\eta)\in M\times\M_p,$  we denote  $f\in C^{E,1}(\M_p)$.
 \item[$(3)$] We denote  $f\in C^{E,1}_K(\M_p)$,  if $f\in C^{E,1}(\M_p)$  and for any compact set $\scr K\subset \M_p$, there exists a constant $C>0$ such that
 $$\sup_{\eta\in \scr K}  | D^E f(\eta)(x)|\le C\big(1+\rr_o^p(x)\big),\ \ x\in M.$$
  \item[$(4)$] We denote  $f\in C^{E,1,1}(\M_p)$, if $f\in C^{E,1}(\M_p)$  such that
 $D^E f(\eta)(x)$ is differentiable in $x$,   $  \nn \{D^Ef(\eta)(\cdot)\}(x)$ is   continuous in $(x,\eta)\in M\times \M_p$, and $|\nn \{D^Ef (\eta)\}|\in L^2(\eta)$ for any $\eta\in \M_p$.
 \item[$(5)$] We write $f\in C^{E,1,1}_B(\M_p)$, if $f\in C^{E,1,1}(\M_p)$ and  for any constant $L>0$ there exists $C_L>0$ such that
  $$\sup_{\eta(\rr_o^p)\le L} |\nn \{D^Ef(\eta)\}|(x)\le C_L(1+\rr_o^p(x)),\ \ x\in M.$$
 \end{enumerate}\end{defn}

Since   for a probability measure $\mu$ and $s>0$, $\mu+s\dd_x$ is no longer a probability measure, for functions of probability measures
we  modify the definition of the extrinsic derivative with the convex combination $(1-s)\mu+s\dd_x$  replacing $\mu+s\dd_x$.

\beg{defn} [{\bf Convexity extrinsic derivative}]    Let $p\in [0,\infty)$ and $f$ be a real function $f$ on $\scr P_p.$
\beg{enumerate} \item[(1)] $f$  is called extrinsically differentiable on $\scr P_p$,  if the centered extrinsic derivative
$$\tt D^Ef(\mu)(x):=\lim_{s\downarrow 0} \ff{f((1-s)\mu+s\dd_x)-f(\mu)}s\in \R $$ exists for all $(x,\mu)\in M\times \scr P_p$.
 \item[$(2)$] We write $f\in C^{E,1}(\scr P_p)$, if   $\tt D^E f(\mu)(x)$ exists and is  continuous in $(x,\mu)\in M\times\scr P_p.$
 \item[$(3)$] We denote  $f\in C^{E,1}_K(\scr P_p)$, if $f\in C^{E,1}(\scr P_p)$  and for any compact set $\scr K\subset \scr P_p$, there exists a constant $C>0$ such that
 $$\sup_{\mu\in \scr K}  | D^E f(\mu)(x)|\le C\big(1+\rr_o^p(x)\big),\ \ x\in M.$$
\item[$(4)$] We write $f\in C^{E,1,1}(\scr P_p)$, if $f\in C^{E,1}(\scr P_p)$ such that $\tt D^E f(\mu)(x)$ is    differentiable in $x\in M$,
  $\nn\{\tt D^Ef(\mu)\}(x)$ is continuous in $(x,\mu)\in M\times\scr P_p$, and $|\nn \{\tt D^Ef(\mu)\}|\in L^2(\mu)$ for any $\mu\in \scr P_p$.
  \item[$(5)$] We write $f\in C_B^{E,1,1}(\scr P_p)$, if $f\in C^{E,1,1}(\scr P_p)$ and for any constant $L>0$ there exists $C>0$ such that
  $$\sup_{\mu(\rr_o^p)\le L} |\nn \{\tt D^Ef(\mu)\}|(x)\le C(1+\rr_o^p(x)),\ \ x\in M.$$
 \end{enumerate} \end{defn}
By Lemma \ref{L2} below with $\gg=\dd_x$ and $r=0$,   we have
$$ \lim_{s\downarrow 0} \ff{f((1-s)\eta+s\dd_x)-f(\eta)}s=  D^Ef(\eta)(x)- \eta\big(D^Ef(\eta)\big),\ \ f\in C_K^{E,1}(\M_p), x\in M.$$
So,   the convexity extrinsic derivative is indeed the centralised extrinsic derivative.

\

 To introduce the intrinsic derivative, for any $v\in \GG_0(TM),$ the class of smooth vector fields  on $M$ with compact support, consider the flow $(\phi_s^v)_{s\ge 0}$ generated by $v$:
 $$\ff{\d}{\d s} \phi_s^v= v(\phi_s^v),\ \ \phi_0^v={\rm Id}, s\ge 0,$$ where {\rm Id} is the identity map.
 Let $\B(TM)$ be the set of all measurable vector fields on $M$. Then for any $\eta\in \M_p$,
 $$L^2(\B(TM);\eta):=\big\{v\in \B(TM): \eta(|v|^2)<\infty\big\}$$ is a Hilbert space, where $|v|:= \ss{\<v,v\>}.$  Then  $\GG_0(TM)$ is dense in $L^2(\B(TM);\eta)$. When $M=\R^d$, we have $\GG_0(TM)= C_0^\infty(\R^d\to\R^d)$ and $\B(TM)=\B(\R^d\to\R^d)$.

 By the Riesz representation theorem,
 for any bounded linear functional $U: \GG_0(TM)\to\R$,
 there exists a unique element $U^*\in L^2(\B(TM);\eta)$ such that
 $$U(v)=\<v, U^*\>_{L^2(\eta)}:=\int_M \<v,U^*\>\d\eta,\ \ v\in \GG_0(TM).$$ In this case,   $U(v):= \int_M \<v,U^*\>\d\eta$ for $ v\in L^2(\B(TM);\eta)$ is  the unique  continuous extension   of $U$ on $L^2(\B(TM);\eta)$.

 \beg{defn}[{\bf Intrinsic derivative}]   Let $p\in [0,\infty)$ and $f$ be a real function on $\M_p$.
  \beg{enumerate} \item[$(1)$] $f$  is called intrinsically differentiable on $\M_p$,    if for any $\eta\in \M_p$ and $v\in \GG_0(TM)$,
 $$D^I_v f(\eta) := \lim_{\vv\downarrow 0} \ff{ f(\eta\circ(\phi_\vv^v)^{-1})-f(\eta)}\vv  \in \R  $$ exists and
 is a bounded linear functional of $v\in \GG_0(TM)\subset L^2(\B(TM);\eta).$ In this case, the unique element $D^If(\eta)\in L^2(\B(TM);\eta)$ such that
 $$D^I_vf(\eta) =\<D^If(\eta),v\>_{L^2(\eta)}:=\int_M \<D^If(\eta), v\> \d\eta,\ \ v\in \GG_0(TM)$$
 is called the intrinsic derivative of $f$ at $\eta\in \M_p.$
\item[$(2)$] We denote  $f\in C^{I,1}(\M_p)$, if $f$ is   intrinsically differentiable   on $\M_p$  such that
 $D^I f $ has a version $D^I f(\eta)(x)$      continuous in $(x,\eta)\in M\times \M_p$.
 %\item[$(3)$] We write $f\in C^{I,1}_b(\M_p)$, if $f\in C^{I,1}(\M_p)$ and   $ D^If $ is bounded.
 \end{enumerate}\end{defn}

 We now introduce the $L$-  and linear functional derivatives  following    \cite{Card, CDLL} where $\scr P_2(\R^d)$ is considered.
  Comparing with the definition of intrinsic derivative, to define the $L$-derivative one   replaces the flow $\phi_s^v$ by  the geodesic flow
  $$\phi_{sv}(x):= \exp_x[sv(x)],\ \ s\ge 0, x\in M,$$
 where
  $$\exp_x: T_xM\to M$$ is the
 exponential map, so that for each $u\in T_xM,$  $$[0,\infty)\ni s\mapsto \gg(s):=\exp_x[su]\in M$$ is the unique geodesic starting from $x$ with initial tangent vector $\ff{\d}{\d s} \gg(s)|_{s=0}=u.$
 When $M=\R^d$, we have $\phi_{sv}(x)= x+sv(x).$ By the triangle inequality, we have
\beq\label{TRAN}  \rr_o(\phi_{sv}(x)):=\rr(o,\phi_{sv}(x))\le \rr(o, x)+\rr(x,\exp_x[sv(x)])\le \rr_o(x)+|v(x)|,\ \ s\in [0,1].\end{equation}  So,
when $p\le 2$, $\eta\in \M_p$ implies
$$\eta\circ\phi_v^{-1}\in \M_p,\ \ v\in L^2(\B(TM);\eta).$$
Thus, in the following definition of $L$-derivative, we assume that $p\le 2$. See   also
\cite{GT19} for a different characterization on the $L$-derivative and applications to the Hamilton-Jacobi equations on $\scr P_2(\R^d).$

\beg{defn}[{\bf $L$-derivative}]  Let $p\in [0,2]$ and $f$ be a real function on $\M_p$.
  \beg{enumerate} \item[$(1)$] $f$  is called weakly $L$-differentiable on $\M_p$,    if for any $\eta\in \M_p$ and $v\in L^2(\B(TM);\eta)$,
 $$D^L_v f(\eta) := \lim_{\vv\downarrow 0} \ff{ f(\eta\circ \phi_{\vv v}^{-1})-f(\eta)}\vv \in \R   $$ exists and 
 is a bounded linear functional of $v\in    L^2(\B(TM);\eta).$ In this case, the unique element $D^Lf(\eta)\in L^2(\B(TM);\eta)$ such that
 \beq\label{DD0} D^L_vf(\eta) =\<D^{L}f(\eta),v\>_{L^2(\eta)},\ \ v\in L^2(\B(TM);\eta)\end{equation}
 is called the weak  $L$-derivative of $f$ at $\eta.$
\item[$(2)$] $f$ is called $L$-differentiable on $\M_p$, if $f$ is weakly $L$-differentiable with
$$\lim_{  \|v\|_{L^2(\eta)}\downarrow 0}   \ff{ |f(\eta\circ \phi_{ v}^{-1})-f(\eta)  -D^{L}_vf(\eta)|}{\|v\|_{L^2(\eta)}} =0,\ \ 0\ne\eta\in \M_p.$$
In this case, we call $D^Lf$ the  $L$-derivative of $f$.
\item[$(3)$] We denote  $f\in C^{L,1}(\M_p)$, if $f$ is   $L$-differentiable   on $\M_p$  such that
 $D^L f $ has a version $D^L f(\eta)(x)$      continuous in $(x,\eta)\in M\times \M_p$.  \end{enumerate}\end{defn}

\beg{defn}[{\bf Linear functional derivative}]   Let $p\in [0,\infty)$ and $f$ be a real function on $\M_p$. A measurable function $$M\ni y\mapsto D^Ff(\eta)(y)$$ is called the linear functional derivative of $f$ at $\eta\in \M_p$, if for any constant $L>0$ there exists a constant $C>0$ such that
\beq\label{LFD1} \sup_{\eta(\rr_o^p)\le L} \Big|D^Ff(\eta)(y)\Big|\le C(1+\rr_o^p(y)),\ \ y\in M,\end{equation}
and for any $\eta,\gg\in \M_p$,
\beq\label{LFD2} f(\gg)-f(\eta)= \int_0^1\d r \int_M D^F f(r\gg+(1-r)\eta)(y)(\gg-\eta)(\d y).\end{equation} Since $(1-s)\mu+s\nu\in \scr P_p$ for $s\in [0,1]$ and $\mu,\nu\in \scr P_p$,  the definition  of $D^F$ also applies to   functions on $\scr P_p$.
%We denote $f\in C^{F,1}(\M_p)$ if $D^Ff(\eta)(y)$ exists and is continuous in $(y,\eta)\in M\times\M_p$.
\end{defn}

The remainder of the paper is organized as follows.   In Section 2,  we state the main results of the paper.   Section 3, we present some lemmas which will be  used in   Sections   4    to prove the main results. The main results of the paper have been reported in the survey \cite{WR20}.

\section{Main results}

 \beg{thm}\label{T1.1}   Let $p\in [0,\infty)$. \beg{enumerate}
 \item[$(1)$] If $f$ is $L$-differentiable on $\M_p$, then it is   intrinsic differentiable and
$  D^If=D^Lf. $
\item[$(2)$] If $f\in C^{E,1}_{K}(\M_p)$, then $f$ has linear functional derivative   $D^Ff= D^Ef.$
 \item[$(3)$]   Let $f\in  C^{E,1,1}(\M_p).$ Then $f \in   C^{I,1}(\M_p)$ with
 \beq\label{FF2}  D^If(\eta)(x)=  \nn \{D^Ef(\eta)(\cdot)\}(x),\ \    (x,\eta)\in M\times\M_p. \end{equation}
 When $p\in [0,2]$ and $f\in  C_B^{E,1,1}(\M_p)$, we have $f\in C^{L,1}(\M_p)$ and
  \beq\label{FF1}  D^Lf(\eta)(x)=  \nn \{D^Ef(\eta)(\cdot)\}(x),\ \    (x,\eta)\in M\times\M_p. \end{equation}
 \item[$(4)$ ] If  $f\in C^{L,1}(\M_p)$, then for any $s\ge 0$,  $ f(\eta+s\dd_\cdot)\in C^1(M)$ with
 \beq\label{AAD} \nn f(\eta+s\dd_{\cdot} )(x)= s D^L f(\eta+s\dd_x)(x),\ \ x\in M, s\ge 0.\end{equation} Consequently,
 \beq\label{FF'} D^L f(\eta)(x)= \lim_{s\downarrow 0} \ff 1 s  \nn f(\eta+s \dd_\cdot)(x),\ \ f\in C^{L,1}(\M_p), (x,\eta)\in M\times\M_p. \end{equation}\end{enumerate}
  \end{thm}

\paragraph{Remark 1.1.}  (a)      Theorem \ref{T1.1}(3) implies   $C_B^{E,1,1}(\M_p)\subset C^{L,1}(\M_p), p\in [0,2].$ However,  a function  $f\in C^{L,1}(\M_p)$ is not necessarily extrinsically differentiable.
For instance, let $\psi\in C([0,\infty))$
but not differentiable, and let  $f(\eta)= \psi(\eta(M))$. Then $f(\eta+s\dd_x)=   \psi(\eta(M)+s)$ which is not differentiable in $s$, so that $f$ is not extrinsically differentiable. But it is easy to see that $f\in C^{L,1}(\M_p)$   with $D^Lf(\eta)=0$. Off course, this counter-example does not work for functions on the space of probability measures

(b) According to  \cite[Proposition 5.48]{CDLL},    if  $f$ is a function on $\scr P_2(\R^d)$ having linear functional derivative $D^Ff(\mu)\in C^1(M)$ for any $\mu\in \scr P_2(\R^d)$,  then $f$ is  $L$-differentiable and
\beq\label{FML} D^{L} f (\eta)= \nn \{D^Ff(\eta)\}.\end{equation}  By Theorem \ref{T1.1}(1)-(3),  this formula \eqref{FML}
is extended to \eqref{FF1}  for the present general framework. Since the definition of $D^E$ is more straightforward than that of $D^F$,   \eqref{FF1}  is more explicit than \eqref{FML}. Note that
in \cite{CDLL} the weak $L$-derivative is named by intrinsic derivative, where the latter   was however  introduced much earlier by \cite{AKR} as in Definition 1.3.

 (c) To illustrate  the link between  derivatives presented in Theorem \ref{T1.1}, let us consider
 the class of cylindrical functions $\F C_b^1$, which consists of functions of type   $$f(\eta):= g(\eta(h_1),\cdots \eta(h_n)),\  \eta(\eta_i):=\int_Mh_i\d\eta,\ \ \ \eta\in\M,$$ where   $n\ge 1, g\in C_b^1(\R^n)$ and $ h_i\in C_b^1(M), 1\le i\le n.$
 Then $f$ is   extrinsically and $L$-differentiable, and has linear functional derivative:
 $$D^L f(\eta)= \sum_{i=1}^n (\pp_i g)(\eta(h_1),\cdots,\eta(h_n)) \nn h_i,\ \
   D^E f(\eta) =D^F f(\eta) = \sum_{i=1}^n (\pp_i g)(\eta(h_1),\cdots,\eta(h_n))   h_i,$$   where  $\nn$ is  the gradient operator on $M$.  Therefore,   we have
  $$ D^Lf(\eta)(x)= \nn \{D^Ef(\eta)(\cdot)\}(x),\ \ (x,\eta)\in M\times\M_p,\ \ p\in [0,\infty) $$
 as indicated in \eqref{FF1}.

  \

  Next, we consider derivatives on the space $\scr P_p:=\M_p\cap\scr P$   for $p\in [0,\infty)$.
Since  for any $\mu\in \scr P_p$ and any $v\in \GG_0(TM)$, we have $\mu\circ\phi_{\vv v}^{-1},\mu\circ(\phi_\vv^v)^{-1}\in \scr P_p$ for $\vv\ge 0$.
 So, the definitions of $D^I$ and $D^L$ work also for functions on $\scr P_p$,
and we define the classes $C^{I,1}(\scr P_p)$ and $C^{L,1}(\scr P_p)$ as in Definitions 1.2 and 1.3 for $\scr P_p$ replacing $\M_p$.

By extending a function on $\scr P_p$ to $\M_p$, we may apply Theorem \ref{T1.1} to establish the corresponding link for functions on $\scr P_2$. As an application,   we will present derivative formula for   the distribution  of random variables.
 For $s_0>0$ and a family of $M$-valued random variables $\{\xi_s\}_{s\in [0,s_0)}$  on a probability space $(\OO,\F,\P),$   we say that $\dot \xi_0:=\ff {\d}{\d s} \xi_s\big|_{s=0}$ exists in
  $L^q(\OO\to TM;\P)$ for some $q\ge 1$,  if  $\dot \xi_0\in T_{\xi_0}M$ with $\E|\dot\xi_0|^q<\infty   $ such that
  \beq\label{*PQ} \lim_{s\downarrow 0} \E\Big| \ff 1 s \exp_{\xi_0}^{-1} [\xi_s] - \dot\xi_0\Big|^q=0.\end{equation}  Since  $\xi_s \to \xi_0$ as $s\to 0$,  note that the inverse of the exponential map
    $\exp_{\xi_0}^{-1} [\xi_s]$   is well-defined for small $s>0$, see the proof of Theorem \ref{T1.1}(1) below for details. In particular, for $M=\R^d$ we have $\exp_{\xi_0}^{-1} [\xi_s]= \xi_s-\xi_0.$

  \beg{cor}\label{T1.2}    Let $p\in [0,\infty)$. \beg{enumerate}
 \item[$(1)$] If $f$ is $L$-differentiable on $\scr P_p$, then it is   intrinsic differentiable and
$D^If=D^Lf. $
\item[$(2)$] If $f\in C^{E,1}_{p}(\scr P_p)$, then $f$ has linear functional derivative on $\scr P_p$ and $D^Ff=\tt D^Ef.$
 \item[$(3)$]  Let $f\in  C^{E,1,1}(\scr P_p).$ Then  $f\in   C^{I,1}(\scr P_p)$ and
 \beq\label{NFF}  D^If(\mu)(x)=   \nn \{\tt D^Ef(\mu)(\cdot)\}(x),\ \    (x,\mu)\in M\times\scr P_p. \end{equation}
 When   $p\le 2$ and $f\in   C_B^{I,1}(\scr P_p)$,  we have  $f\in C^{L,1}(\scr P_p)$ with
 $$D^Lf=\nn \{\tt D^Ef(\mu)(\cdot)\},\ \ \mu\in\scr P_p,f\in C_B^{E,1,1}(\scr P_p). $$
 \item[$(4)$ ] If  $f\in C^{L,1}(\scr P_p)$, then  $ f((1-s)\mu+s\dd_\cdot)\in C^1(M)$ with
 \beq\label{NAAD} \nn f((1-s)\mu+s\dd_{\cdot} )(x)= s D^L f((1-s)\mu+s\dd_x)(x),\ \ x\in M.\end{equation} Consequently,
 \beq\label{NFF'} D^L f(\mu)(x)= \lim_{s\downarrow 0} \ff 1 s  \nn f((1-s)\mu+s \dd_\cdot)(x),\ \ f\in C^{L,1}(\scr P_p), (x,\eta)\in M\times\M. \end{equation}
 \item[$(5)$] Let  $\{\xi_s\}_{s\in [0,s_0)}$  be random variables on $M$ with $ \L_{\xi_s}  \in \scr P_p$  continuous in $s$,  such that $\dot \xi_0:=\ff {\d}{\d s} \xi_s\big|_{s=0}$ exists in
  $L^q(\OO\to TM;\P)$ for some $q\ge 1$. Then
    \beq\label{LLP0}  \lim_{s\downarrow 0}\ff{f(\L_{\xi_s})-f(\L_{\xi_0})} s = \E \big\<D^Lf(\L_{\xi_0}) (\xi_0),\dot \xi_0\big\>\end{equation}
    holds for any $f\in C^{E,1,1}(\scr P_p)$ such that for any compact set $\scr K\subset \scr P_p$,
   \beq\label{*10} \sup_{\mu\in \scr K} |\nn \{\tt D^E f(\mu)\}|(x)\le C(1+\rr_o)^{\ff{p(q-1)}q},\ \ x\in M\end{equation} holds for some constant $C>0$.
 \end{enumerate}
  \end{cor}

Let us  compare \eqref{LLP0}  with the corresponding  formula    presented in  \cite{Card}  for $M=\R^d,\rr(x,y)=|x-y|$ and $p=2$.
In this case,    the formula \eqref{LLP0} is established   for the  probability space being Polish and $f\in C^{L,1}(\scr P_2(\R^d))$ with bounded $D^Lf$,
see also \cite[Proposition A.2]{HSS} and \cite[Lemma 2.3]{W19} for this formula with more  general functions $f$ on $\scr P_2(\R^d)$.    Theorem  \ref{T1.2}   establishes  \eqref{LLP0} to $\M_p$ on  Riemannian manifolds and $p\ge 0$.

  \section{Some lemmas}

We first consider the variation of $f(h\eta)$ in the density function $0\le h\in L^1(\eta)$. Recall that for a nonnegative measurable function $h$ on $M$,   the measure $h \eta $ is defined by
 $$(h\eta)(A):= \int_A h \d\eta,\ \ A\in \B(M),$$ where $\B(M)$ is the Borel $\si$-field of $M$.
  In the proof of Theorem \ref{T1.1}, we will formulate $f(\eta\circ \phi_{\vv v}^{-1})$ by
  $  f((1+ h_\vv) \eta )$  for some $h\in \scr H_{\vv_0}$, where $h\in \scr H_{\vv_0}$ means that $h\in C([0,\vv_0]\times M; [0,\infty))$ and
 \beg{enumerate} \item[(1)] $h_0=0$,  $\sup_{\vv\in [0,\vv_0]}\|h_\vv\|_\infty<\infty,$  $\supp h_\vv\subset K$ for some compact set $K\subset M$ and all $\vv\in [0,\vv_0]$;
 \item[(2)] $\dot h_\vv:=\lim_{s\downarrow 0} \ff{h_{\vv+s}-h_\vv}s\in C_b(M)$ exists and is uniformly bounded for $\vv\in [0,\vv_0)$.    \end{enumerate}
So, to calculate $D^Lf(\mu)$,   we  first present the following lemma which  links $ f((1+h_\vv)\eta)- f(\eta)$  to the extrinsic derivative.

 \beg{lem} \label{L1} Let $p\in [0,\infty)$.  For any $ h\in \scr H_{\vv_0}$  and any $f\in C^{E,1,1}(\M_p)$,
  \beq\label{DF}   f((1+h_\vv)\eta) -f(\eta) = \int_0^\vv\d r  \int_{M} D^Ef((1+h_r)\eta)(x) \dot h_r(x) \eta(\d x),\ \ \eta \in \M_p,\vv\in [0,\vv_0].\end{equation}
 \end{lem}

 \beg{proof} (1) We first consider   $$\eta\in \M_{disc}
   :=\Big\{\sum_{i=1}^n a_i\dd_{x_i}: n\ge 1, a_i>0, x_i\in M, 1\le i\le n\Big\}. $$
 In this case, for   any $\vv\in [0,\vv_0)$ and $s\in (0,\vv_0-\vv)$, by the definition of $D^E$ we have
 \beg{align*} &f((1+  h_{\vv+s}) \eta ) - f((1+h_\vv)\eta)= f\Big((1+h_{\vv})\eta+ \sum_{i=1}^n \{h_{\vv +s}-h_\vv\}(x_i)a_i \dd_{x_i}\Big) -f((1+h_\vv)\eta)\\
 &= \sum_{k=1}^n \bigg\{f\Big((1+h_\vv)\eta+ \sum_{i=1}^k \{h_{\vv +s}-h_\vv\}(x_i)a_i \dd_{x_i}\Big) -f\Big((1+h_\vv)\eta+ \sum_{i=1}^{k-1} \{h_{\vv +s}-h_\vv\}(x_i)a_i \dd_{x_i}\Big)\bigg\}\\
 &= \sum_{k=1}^na_k \int_{-a_k   \{h_{\vv +s}-h_\vv\}^-(x_k)}^{a_k   \{h_{\vv +s}-h_\vv\}^+(x_k)} \Big\{D^E f\Big((1+h_\vv)\eta+ \sum_{i=1}^{k-1} \{h_{\vv +s}-h_\vv\}(x_i)a_i \dd_{x_i}+  r \dd_{x_k}\Big)\Big\}(x_k)\d r,\end{align*}
 where  $\sum_{i=1}^0:=0$,
  $a^+:= \max\{a,0\}$ and $ a^-:= (-a)^+$ for $a\in\R$. Multiplying by $s^{-1}$ and letting $s\downarrow 0$, we deduce from this and the continuity of $D^Ef$ that
\beq\label{PL} \beg{split} & \lim_{s\downarrow 0} \ff {f((1+  h_{\vv+s}) \eta ) - f((1+h_\vv)\eta) }s   = \sum_{k=1}^n a_k \{\dot h_\vv(x_k)^+-\dot h_\vv(x_k)^-\}D^Ef((1+h_\vv)\eta)(x_k)\\
&=\int_{M} D^Ef((1+h_\vv)\eta)(x) \dot h_\vv(x) \eta(\d x),\ \ \vv\in [0,\vv_0), \eta\in\M_{disc}.\end{split}\end{equation}

(2) In general, for any $\eta\in \M_p$, let $\{\eta_n\}_{n\ge 1} \subset \M_{disc}$ such that $\eta_n\to\eta$ in $\M_p$.  By \eqref{PL}, for any $\vv\in (0,\vv_0)$ and $s\in (0,\vv_0-\vv)$,    we have
\beq\label{NM}  f((1+  h_{\vv}) \eta_n)-f(\eta_n)  =    \int_0^\vv\d r \int_M  D^E f ((1+h_r) \eta_n)(x) \dot h_r(x) \eta_n(\d x),\ \ n\ge 1.\end{equation}
Next, since  $D^Ef\in C(M\times\M_p)$ and $h_r,\dot h_r\in C_b(M)$ for $r\in [0, \vv_0]$ with compact support $\subset K$,  and $\eta_n\to\eta$ in $\M_p$, we obtain
\beq\label{NM2} \lim_{n\to\infty} \int_M D^Ef((1+h_r)\eta )(x) \dot h_r(x) \eta_n(\d x) = \int_M D^Ef((1+h_r)\eta)(x) \dot h_r(x) \eta(\d x).\end{equation}
Moreover,   $\eta_n\to \eta$ in $\M_p$ and $h\in \scr H_{\vv_0}$ imply that the set $$\scr K_r:=\{(1+h_r)\eta, (1+h_r)\eta_n:n\ge 1 \}$$  is compact in $\M_p$ for any $r\in [0,\vv_0]$. Combining this with   $D^Ef\in C(M\times\M_p)$, we see that    the function
$$\scr K_r\times M\ni (\gg,x)\mapsto D^Ef(\gg )(x) \dot h_r(x)$$ is uniformly continuous and has compact support $\subset \scr K_r\times K$, so that   \eqref{NM2} implies
\beg{align*}&\limsup_{n\to \infty} \bigg|\int_M D^Ef((1+h_r)\eta_n )(x) \dot h_r(x) \eta_n(\d x) - \int_M D^Ef((1+h_r)\eta)(x)\dot h_r(x) \eta(\d x)\bigg|\\
  &= \limsup_{n\to \infty} \bigg|\int_M D^Ef((1+h_r)\eta_n )(x)\dot h_r(x) \eta_n(\d x) - \int_M D^Ef((1+h_r)\eta)(x)\dot h_r(x) \eta_n(\d x)\bigg|\\
  &\le   \limsup_{n\to \infty} \Big\{\eta_n(K)  \sup_{x\in K} |D^Ef((1+h_r)\eta_n )(x)\dot h_r(x)  -   D^Ef((1+h_r)\eta)(x)\dot h_r(x) |  \Big\}\\
  &=0.\end{align*}
  Combining this with
  $$\sup_{(\gg,x)\in \scr K_r\times K,r\in [0,\vv_0]} |D^E f(\gg)(x) \dot h_r(x)|<\infty,$$
  we deduce from   the dominated convergence theorem that
 \beq\label{**0} \beg{split} &\lim_{n\to\infty} \int_0^\vv\d r  \int_M \big\{D^E f\big\}((1+h_\vv) \eta_n)(x) \dot h_r(x) \eta_n(\d x)\\
 &=  \int_0^\vv\d r \int_M \big\{D^E f\big\}((1+h_r) \eta)(x) \dot h_r(x) \eta(\d x).\end{split} \end{equation}
 Therefore, by letting $n\to\infty$ in \eqref{NM} and using the continuity of $f$, we prove \eqref{DF}.
\end{proof}

To calculate the convexity extrinsic derivative, we present the following result.

\beg{lem}\label{L2} Let $p\in [0,\infty)$. Then for any  $f\in C^{E,1}_K(\M_p)$ and  $\eta,\gg\in \M_p$,
\beg{align*} &\ff{\d }{\d r} f((1-r)\eta+ r\gg) := \lim_{\vv\downarrow 0} \ff{f((1-r-\vv)\eta+ (r+\vv)\gg) -f((1-r)\eta+ r\gg) }\vv \\
&=\int_M \big\{D^Ef((1-r)\eta+ r\gg) (x)\big\}(\gg-\eta)(\d x),\ \ r\in [0,1).\end{align*} Consequently, for any $f\in C^{E,1}_K(\M_p)$,
\beg{align*} &\tt D^Ef(\eta)(x):= \lim_{s\downarrow 0} \ff{f((1-s)\eta+s\dd_x)-f(\eta)}s\\
 &= D^Ef(\eta)(x)-\eta\big(D^Ef(\eta)\big),\ \   (x,\eta)\in M\times \M_p.\end{align*}
The assertions also hold for $\scr P_p$ replacing $\M_p$. \end{lem}

 \beg{proof}  As in the proof of Lemma \ref{L1}, we  take
  $$\eta_n=\sum_{i=1}^n\aa_{n,i}\dd_{x_{n,i}}, \ \ \gg_n=\sum_{i=1}^n \bb_{n,i}\dd_{x_{n,i}}$$ for some $x_{n,i}\in M$ and $\aa_{n,i},\bb_{n,i}\ge 0$, such that
 $$\eta_n\to \eta,\ \ \gg_n\to\gg  \ \ \text{ in}\  \M_p\ \text{as}\  n\to\infty.$$
 For any $r\in [0,1)$ and $\vv\in (0,1-r)$,  let
 $$\LL_{n,i}^\vv:= (1-r)\eta_n+ r\gg_n + \sum_{k=1}^{i-1} \vv(\bb_k-\aa_k)\dd_{x_{n,k}}\in\M_p,\ \ 1\le i\le n,$$ where by convention $\sum_{i=1}^0:=0$. Then by the definition of $D^Ef$, we have
 \beg{align*} &f((1-r-\vv)\eta_n+ (r+\vv)\gg_n) -f((1-r)\eta_n+ r\gg_n) \\
 &=\sum_{i=1}^n  \big\{f(\LL^\vv_{n,i}+ \vv(\bb_{n,i}-\aa_{n,i})\dd_{x_{n,i}}) -f(\LL^\vv_{n,i})\big\} \\
 &= \sum_{i=1}^n \int_{-\vv(\bb_{n,i}-\aa_{n,i})^- }^{\vv(\bb_{n,i}-\aa_{n,i})^+}  D^E f(\LL^\vv_{n,i}+ s\dd_{x_{n,i}})  (x_{n,i}) \d s,\ \ \vv\in (0,1-r).\end{align*}
Multiplying by $\vv^{-1}$ and letting $\vv\downarrow 0$, due to the continuity of $D^Ef$ we derive
\beg{align*} &\ff{\d }{\d r} f((1-r)\eta_n+ r\gg_n)= \sum_{i=1}^n (\bb_{n,i}-\aa_{n,i})  D^Ef((1-r)\eta_n+ r\gg_n) (x_{n,i})  \\
&=\int_M \big\{D^Ef((1-r)\eta_n+ r\gg_n) (x)\big\}(\gg_n-\eta_n)(\d x),\ \ r\in [0,1),\ \ n\ge 1.\end{align*}
Consequently, for any $r\in [0,1)$,
\beg{align*} &f((1-r-\vv)\eta_n+ (r+\vv)\gg_n) -f((1-r)\eta_n+ r\gg_n)\\
&=\int_{r}^{r+\vv} \d s \int_M \big\{D^Ef((1-s)\eta_n+ s\gg_n) (x)\big\}(\gg_n-\eta_n)(\d x),\ \ \vv\in (0,1-r), n\ge 1.\end{align*}
Noting that the set $\{\eta_n,\gg_n: n\ge 1\}$ is relatively compact in $\M_p$, by this and the condition on $f$, we may let $n\to\infty$ to derive
\beg{align*} &f((1-r-\vv)\eta+ (r+\vv)\gg) -f((1-r)\eta+ r\gg)\\
&=\int_{r}^{r+\vv} \d s \int_M \big\{D^Ef((1-s)\eta+ s\gg) (x)\big\}(\gg-\eta)(\d x),\ \ \vv\in (0,1-r).\end{align*}
Multiplying by $\vv^{-1}$ and letting $\vv\downarrow 0$, we finish the proof. \end{proof}

The following is a consequence  of Lemma \ref{L2} for functions on $\scr P_p$.
\beg{lem}\label{LL'} Let $p\in [0,\infty)$. Then for any  $f\in C^{E,1}_K(\scr P_p)$ and  $\mu,\nu\in \scr P_p$,
\beg{align*}   \lim_{s\downarrow 0} \ff{f((1-s)\mu+ s\nu) -f(\mu) }s
 =\int_M \big\{\tt D^Ef((\mu) (x)\big\}(\nu-\mu)(\d x).\end{align*}   \end{lem}

\beg{proof}
To apply Lemma \ref{L2},  we extend a function $f$ on $\scr P_p$  to $\tt f$ on $\M_p$   by letting
 $$\tt f(\eta)= h(\eta(M)) f(\eta/\eta(M)),\ \ \eta\in \M_p,$$
 where $h\in C^\infty_0(\R)$ with support contained by $[\ff 1 4,2]$ and $h(r)=1$ for $r\in [\ff 1 2,\ff 3 2].$ It is easy to see that
 $$f((1-s)\mu+s\nu)= \tt f((1-s)\mu+s\nu),\ \ s\in [0,1], \mu,\nu\in \scr P_p,$$
 and  $f\in C_K^{E,1}(\scr P_p)$ implies that $\tt f\in C_K^{E,1}(\M_p)$ and
 $$ D^E\tt f(\mu)= \tt D^Ef(\mu),\ \   \mu\in \scr P.$$ Then the desired formula is implied by   Lemma \ref{L2} with $r=0$.\end{proof}

Finally, we prove a derivative formula for the distribution of random variables.

 \beg{lem}\label{L3}   Let  $\{\xi_s\}_{s\in [0,s_0)}$  be  $M$-valued random variables  such that $ \lim_{s\to 0} \L_{\xi_s}=\L_{\xi_0}$ in $\scr P_p$, and $\dot \xi_0:=\ff {\d}{\d s} \xi_s\big|_{s=0}$ exists in
  $L^q(\OO\to TM;\P)$ for some $q\ge 1$. Then
   \beq\label{LLP}  \lim_{s\downarrow 0} \ff{f(\L_{\xi_s})-f(\L_{\xi_0})} s = \E \big\<\nn\{\tt D^E f(\L_{\xi_0}) \}(\xi_0),\dot \xi_0\big\> \end{equation}    holds for functions  $f\in C^{E,1,1}(\scr P_p)$ satisfying \eqref{*10}
    for any compact $\scr K\subset \scr P_p$ and some constant $C=C(\scr K)>0.$
    \end{lem}
 \beg{proof} By Lemma \ref{LL'}, we have
 \beq\label{AV} \beg{split}  & f(\L_{\xi_{s}})-f(\L_{\xi_0})= \int_0^1\Big\{ \ff{\d}{\d r} f(r \L_{\xi_{s}}+ (1-r) \L_{\xi_0})\Big\}\d r \\
 &= \int_0^1\d r \int_M \big\{\tt D^Ef (r \L_{\xi_{s}}+ (1-r) \L_{\xi_0})\big\}(x) (\L_{\xi_s}-\L_{\xi_0})(\d x) \\
 &= \int_0^1\E\Big[\big\{\tt D^E f (r \L_{\xi_{s}}+ (1-r) \L_{\xi_0})\big\}(\xi_s) -\big\{\tt D^E f (r \L_{\xi_{s}}+ (1-r) \L_{\xi_0})\big\}(\xi_0) \Big]\d r. \end{split}\end{equation}
 For each $s\ge 0,$ let $\gg_{s,\cdot}: [0,1]\to M$ be the minimal geodesic such that
 $\gg_{s,0}=\xi_0$ and $\gg_{s,1}=\xi_s.$ Then $\lim_{s\downarrow  0} \gg_{s,\theta}= \xi_0$, and   by \eqref{*PQ},
 $$ \lim_{s\downarrow 0}  \E\bigg|\ff 1 s  //_{\theta\to 0}\ff{\d}{\d\theta} \gg_{s,\theta}-\dot \xi_0\bigg|^q=0, $$
where $//_{\theta\to 0}: T_{\gg_{s,\theta}}M\to T_{\xi_0}M$ is the parallel displacement along the geodesic $\gg_{s,\cdot}: [0,\theta]\to M$.
Combining these with  \eqref{AV} and \eqref{*10} with $\scr K:=\{\L_{\xi_0}, \L_{\xi_{s_n}}: n\ge 1\}$ for a sequence $s_0> s_n\downarrow 0$,   we may apply the dominated convergence theorem to derive
\beg{align*} & \lim_{s_n\downarrow 0}  \ff{f(\L_{\xi_{s_n}})-f(\L_{\xi_0})} {s_n} \\
&=  \lim_{s_n\downarrow 0} \ff 1 {s_n} \int_0^1\E\Big[\big\{D^E f (r \L_{\xi_{s_n}}+ (1-r) \L_{\xi_0})\big\}(\xi_{s_n}) -\big\{\tt D^E f (r \L_{\xi_{s_n}}+ (1-r) \L_{\xi_0})\big\}(\xi_0) \Big]\d r\\
&= \lim_{s_n\downarrow 0}  \int_0^1\d r \int_0^1 \E\Big[\big\<\nn\big\{\tt D^E f(r \L_{\xi_{s_n}}+ (1-r) \L_{\xi_0})\big\}(\gg_{s_n,\theta}),
\ff 1 s  \ff{\d}{\d\theta} \gg_{s_n,\theta}\big\>\Big]\d\theta\\
&= \E \big\<\nn\{\tt D^E f(\L_{\xi_0}) \}(\xi_0),\dot \xi_0\big\>.\end{align*}
\end{proof}

 \section{Proofs of Theorem \ref{T1.1} and Corollary \ref{T1.2}  }
Obviously, assertion (2) follows from Lemma \ref{L2}. Below we prove assertions (1),  (3),  (4) in Theorem \ref{T1.1} as well as Corollary \ref{T1.2} respectively.

\beg{proof}[Proof of Theorem $\ref{T1.1}(1)$] Although the flows $\phi_s^v$ and $\phi_{sv}$ are different, their derivative at $s=0$ are all equal to $v$, so that both $D^I_v$ and $D^L_v$ are directional derivatives along $v$. Thus, it is reasonable that for a large class of functions we have $D^If=D^Lf$. To see this, we need the inverse exponential map $\exp_x^{-1}.$  For any $z\in M$, let $u\in T_xM$ such that
 $$[0,1]\ni s\mapsto \exp_x[s u]\in M$$ is the minimal geodesic from $x$ to $z$, and we denote $u=\exp_x^{-1}[z].$ If $z$ is not in the cut-locus of $x$, the minimal geodesic from $x$ to $z$ is unique, and $\exp_x^{-1}[z]$
 is smooth in $z$.
 In case that $z$ belongs to the cut-locus of $x$,
 such a vector $u\in T_xM$ may be not unique. For any compact set $\scr K\subset M$, there exists a constant $R >0$ such that
 for any $x\in \scr K$, the distance between $x$ and its cut-locus is larger than $R $. So, for any $x\in \scr K$,
 $$\exp_x: \{u\in T_xM: |u|\le R\}\to B_x(R ):= \{y\in M: \rr(x,y)\le R\}$$ is a diffeomorphism, such that
 $$\exp_x^{-1}: B_x(R)\to T_xM$$ is smooth.
Thus,   for any $v\in \GG_0(TM)$ and small enough $\vv>0$,   we have $v_\vv:= \exp_x^{-1}[\phi_\vv^v] \in \GG_0(TM)$. Moreover,
$$v_\vv=   \vv v+{\rm o}(\vv),$$
where $ \vv^{-1}\|{\rm o}(\vv)\|_\infty\to 0$ as $\vv\downarrow 0$. Hence, for any $L$-differentiable function $f$ and  $\eta\in \M$, when $\vv$ is small enough we have
\beg{align*} &\limsup_{\vv\downarrow 0} \bigg|\ff{f(\eta\circ(\phi_\vv^v)^{-1})-f(\eta)}\vv- D^L_vf(\eta)\bigg|
= \limsup_{\vv\downarrow 0}\bigg|\ff{f(\eta\circ \phi_{v_\vv}^{-1})-f(\eta)}\vv- D^L_vf(\eta)\bigg|\\
&\le \limsup_{\vv\downarrow 0}   \bigg|\ff{f(\eta\circ \phi_{v_\vv}^{-1})-f(\eta) - D^L_{v_\vv} f(\eta)}{\vv}\bigg|+ \big|D^L_{v-\vv^{-1}  v_\vv}f(\eta)\big|\bigg\}=0.\end{align*}
Therefore, $D^If=D^Lf$ holds for $L$-differentiable $f$. \end{proof}

  \beg{proof}[Proof of Theorem $\ref{T1.1}(3)$]  It suffices to prove the formulas \eqref{FF2} and \eqref{FF1} for $f\in C^{E,1,1}(\M_p)$
  and $f\in C^{E,1,1}(\M_p)$ respectively.

(a) For \eqref{FF2}. Since     any  $\eta\in \M_p$   can be approximated  by those having smooth and strictly positive density functions with respect to the volume measure $\d x$,
  by the argument leading to \eqref{**0}, it suffices to show that for any $\eta\in \M_p$ satisfying
 \beq\label{ETA} \eta(\d x) = \rr(x)\d x\ \text{for\ some}\  \rr\in C^\infty_b(M),\ \inf \rr>0,\end{equation}   there exists  a constant $\vv_0>0$ such that
    \beq\label{DG'}     f(\eta\circ (\phi_{\vv}^v)^{-1}) -f(\eta) =  \int_0^\vv \d r  \int_M \<\nn\{D^Ef(\eta\circ\phi_{r v}^{-1})\}, v\>\d(\eta\circ\phi_{r v}^{-1}),\ \ \vv\in (0,\vv_0).\end{equation}
Firstly, there exists a constant $\vv_0>0$ such that
  $$\rr_{\vv}^v:=\ff{\d(\eta\circ(\phi_{\vv}^v)^{-1}) }{\d \eta },\ \  \dot \rr_{\vv}^v:= \lim_{s\downarrow 0}  \ff{ \rr_{\vv+s}^v-\rr_{\vv}^v}{s} $$ exist in $C_b(M)$ and are uniformly bounded and continuous in $\vv\in [0,\vv_0]$.
Next,  by Lemma \ref{L1},  we have
 \beq\label{NN'}    f(\eta\circ (\phi_{\vv}^v)^{-1})-f(\eta)
 = \int_0^\vv\d r \int_{M} \big\{D^E f(\eta\circ(\phi_{r}^v)^{-1})\big\}   \dot \rr_{r}^v\d\eta,\ \ \vv\in [0,\vv_0]. \end{equation}
To calculate $\dot\rr^v_{r}$, by $\ff{\d}{\d s} \phi_s^v= v(\phi_s^v)$,   for any $g\in C_0^\infty(M)$ we have
 $$\ff{\d}{\d r} \big\{g\circ\phi_{r}^v\big\}= \<\nn g (\phi_{r}^v),  v(\phi_{r}^v)\>= \<\nn g, v\>(\phi_{r}^v),\ \ \ r\ge 0,$$ which is smooth and   bounded in $(r,x)\in [0,\vv_0]\times M$. So,
\beg{align*} &\int_M g \dot \rr_r^v \d\eta= \int_M g \lim_{s\downarrow 0} \ff{\rr^v_{r+s}- \rr^v_r} s \d\eta =\lim_{s\downarrow 0} \ff 1 s \int_M g \d\big\{\eta\circ(\phi_{r+s}^v)^{-1}- \eta\circ(\phi_{r}^v)^{-1}\big\}\\
&=\lim_{s\downarrow 0} \ff 1 s \int_M \big\{g\circ \phi_{r+s}^v - g\circ\phi_{r}^v  \big\}\d\eta = \int_M \ff{\d}{\d r} (g\circ \phi_r^v)\,\d\eta\\
&=\int_M \<\nn g,v\>\circ\phi_{r} ^v\d \eta= \int_M \<\nn g,v\>  \d (\eta\circ(\phi_{\vv}^v)^{-1}) \\
&=- \int_M \big\{g \,\div_{\eta\circ(\phi_{r}^v)^{-1}} (v) \big\}  \d (\eta\circ(\phi_{r}^v)^{-1})= - \int_Mg \big\{\div_{\eta\circ(\phi_{r}^v)^{-1}} (v)  \rr_r^v\big\}\d  \eta,\ \ g\in C_0^\infty(M),\end{align*} where $\div_{\eta\circ(\phi_{r}^v)^{-1}} (v)= \div (v)+ \<v,\nn\log (\rr_r^v\rr)\>.$
This implies  $\dot \rr_r^v= -\div_{\eta\circ(\phi_{r}^v)^{-1}} (v)  \rr_r,$ so that the integration by parts formula and $\rr_r^v\eta= \eta\circ(\phi_r^v)^{-1}$ lead to
\beg{align*} &\int_M \big\{D^E f(\eta\circ(\phi_{r}^v)^{-1})\big\}\dot\rr_r^v\d\eta = - \int_M \big\{D^E f(\eta\circ(\phi_{r}^v)^{-1})\big\}{\rm div}_{\eta\circ(\phi_r^v)^{-1}}(v)\, \d (\eta\circ(\phi_r^v)^{-1}) \\
&= \int_M\big\<\nn \{D^E f(\eta\circ(\phi_{r}^v)^{-1})\},v\big\>\d (\eta\circ(\phi_{r}^v)^{-1}).\end{align*}
Combining this with   \eqref{NN'} we prove  \eqref{DG'}.

(b) For \eqref{FF1}. Let $p\in [0,2].$ For any $\eta\in \M_p$ and
$v\in L^2(\B(TM);\eta)$ with $\eta(|v|^2)\le 1$, by \eqref{TRAN}  we have
$$\sup_{s\in [0,1]} (\eta\circ\phi_{sv}^{-1})(\rr_o^p) =\eta(\rr_o(\phi_{sv})^p)\le 2 \eta(\rr_o^p+ |v|^p)<\infty.$$
Then  there exists a constant $K>0$ such that
\beq\label{DDN} \sup_{s\in [0,1],\eta(|v|^2)\le 1} (\eta\circ\phi_{sv}^{-1}+\eta)(\rr_o^p)\le K.\end{equation}
So, by Lemma \ref{L2}, we obtain
 \beg{align*} &  f(\eta\circ\phi_v^{-1})- f(\eta)= \int_0^1 \Big\{\ff{\d}{\d r} f(r \eta\circ\phi_{rv}^{-1} +(1-r)\eta)\Big\}\d r\\
 &= \int_0^1 \d r\int_M  (D^Ef)(r \eta\circ\phi_{rv}^{-1} +(1-r)\eta) \d(\eta\circ\phi_{v}^{-1}-\eta) \\
 &= \int_0^1 \d r\int_M    \Big\{(D^Ef)(r \eta\circ\phi_{rv}^{-1} +(1-r)\eta)(\phi_{v}(x))-(D^Ef)(r \eta\circ\phi_{rv}^{-1} +(1-r)\eta)(x) \Big\} \eta(\d x)\\
  &= \int_0^1 \d r\int_M  \eta(\d x) \int_0^1 \big\<//_{\phi_{sv}(x)\to x}\nn \big\{(D^Ef)(r \eta\circ\phi_{rv}^{-1} +(1-r)\eta)\big\}(\phi_{sv}(x)),v(x)\big\> \d s,\end{align*}
  where $//_{\phi_{sv}(x)\to x}: T_{\phi_{sv}(x)}M\to T_xM$ is the parallel displacement along the geodesic $[0,s]\ni\theta\mapsto \phi_{(s-\theta)v}(x).$ Thus,
  \beg{align*} &I_v:= \ff{|f(\eta\circ\phi_v^{-1})- f(\eta)-\int_M\<\nn\{D^Ef(\eta)\},v\>\d\eta|^2}{\eta(|v|^2)}\\
  &\le  \int_{[0,1]^2\times M}    \big|//_{\phi_{sv}(x)\to x}\nn \big\{(D^Ef)(r \eta\circ\phi_{rv}^{-1} +(1-r)\eta)\big\}(\phi_{sv}(x))-\nn\big\{D^Ef(\eta)\big\}(x)\big|^2\d r\d s\eta(\d x).\end{align*}
  By \eqref{DDN}, as $\|v\|_{L^2(\eta)}\to 0$ we have $\phi_{sv}(x)\to x\ \eta$-a.e. and
  $\eta\circ\phi_{sv}^{-1}\to \eta$ in $\M_p$ for any $s\ge 0$. Combining these with \eqref{DDN} we may apply the dominated convergence theorem to derive $I_v\to 0$ as $\|v\|_{L^2(\eta)}\to 0$. Therefore, $f$ is $L$-differentiable such that \eqref{FF1} holds.

\end{proof}

 \beg{proof}[Proof of Theorem $\ref{T1.1}(4)$]   It suffices to prove \eqref{AAD}. Let  $f\in C^{L,1}(\M)$.
 We first prove the formula for $\eta\in \M_p$ and $x\in M$ with $\eta(\{x\})=0$, then extend to the general situation.

(a) Let  $\eta(\{x\})=0$.  In this case, for any $v_0\in T_xM$,   let $v=1_{\{x\}}v_0.$ Then
$$\phi_{rv}(z)= \beg{cases} z,\ &\text{if}\ z\ne x,\\
\exp_x[rv_0],\ &\text{if}\ z= x.\end{cases}$$  By $\eta(\{x\})=0$, we have
\beq\label{PLL} (\eta+s\dd_x)\circ \phi_{rv}^{-1}= \eta+ s\dd_{\exp_x[rv_0]}.\end{equation}Since   $v$ can be approximated in $L^2(\eta+s\dd_x)$ by elements in $\GG_0(TM)$, the $L$-differentiability of $f$ and $\eta(\{x\})=0$ imply
\beg{align*}&\lim_{r\downarrow 0} \ff{f((\eta+s\dd_x)\circ\phi_{rv}^{-1})-f(\eta+s\dd_x)}r \\
&=\int_M  \<D^Lf(\eta+s\dd_x),v\> \d (\eta+s\dd_x))= s\<D^Lf(\eta+s\dd_x)(x),v_0\>.\end{align*}
Combining this with \eqref{PLL}, we obtain \beg{align*}   \lim_{r\downarrow 0} \ff{ f(\eta+s\dd_{\exp_x[rv_0]})- f(\eta+s\dd_x)}r= s \<D^Lf(\eta+s\dd_x)(x), v_0\>.\end{align*}
This implies that  $f(\eta+s\dd_\cdot)$ is differentiable at point $x$ and \eqref{AAD} holds.

(b) In general, for any $v_0\in T_xM$, there exists $r_0>0$ such that   $v_0$ extends to a smooth vector field  $v$ on $B(x,r_0)$    by parallel displacement; i.e. $v(x)$ is the parallel displacement along the minimal geodesic from $x$ to $z$.
Since $\eta(\{\exp_x[\theta v_0]\})=0$ for a.e. $\theta \ge 0$, by the continuity of $f$ and the formula \eqref{AAD} for $\eta(\{x\})=0$ proved above, we obtain
\beg{align*}  & \ff{ f(\eta+s\dd_{\exp_x[rv_0]})- f(\eta+s\dd_x)}r = \ff 1 r \int_0^r \ff{\d}{\d\theta} f(\eta+s\dd_{\exp_x[\theta v_0]}) \d\theta\\
&= \ff 1 r \int_0^r\big\< \nn  f(\eta+s\dd_{\cdot}) (\exp_x[\theta v_0]), v\big(\exp_x[\theta v_0]\big)\big\> \d\theta\\
&=  \ff s r \int_0^r \big\<D^L f(\eta+s\dd_{\cdot}) (\exp_x[\theta v_0]), v\big(\exp_x[\theta v_0]\big)\big\>  \d\theta,\ \ r\in (0,r_0).\end{align*}
By the continuity of $D^Lf$, with $r\downarrow 0$ this implies \eqref{AAD}.\end{proof}

 \beg{proof}[Proof of Corollary  \ref{T1.2}]   To apply Theorem \ref{T1.1},  we extend a function $f$ on $\scr P_p$  to $\tt f$ on $\M_p$ as in the proof of Lemma \ref{LL'}, i.e. by letting
 $$\tt f(\eta)= h(\eta(M)) f(\eta/\eta(M)),\ \ \eta\in \M_p,$$
 where $h\in C^\infty_0(\R)$ with support contained in $[\ff 1 4,2]$ and $h(r)=1$ for $r\in [\ff 1 2,\ff 3 2].$ It is easy to see that
 $$f((1-s)\mu+s\nu)= \tt f((1-s)\mu+s\nu),\ \ s\in [0,1], \mu,\nu\in \scr P_p,$$
 and  $f\in C^{E,1,1}(\scr P_p)$ implies that $\tt f\in C^{E,1,1}(\M_p)$ and
 $$D^E\tt f(\mu)= \tt D^Ef(\mu),\ \ D^Lf(\mu)= D^L\tt f(\mu),\ \ D^If(\mu)=D^{int } \tt f(\mu), \ \  \mu\in \scr P.$$
Then Corollary \ref{T1.2}(1)-(4) follow from the corresponding assertions in Theorem \ref{T1.1}  with $\tt f$ replacing $f$.

Finally,   since
 $\nn\{\tt D^Ef(\mu)\}= \nn\{D^E \tt f\}(\mu)= D^Lf(\mu)$ for $\mu\in \scr P_p$ and $f\in C^{E,1,1}(\scr P_p)$,      \eqref{LLP0}  follows from Lemma \ref{L3}.
   \end{proof}

\paragraph{Acknowledgement.} The authors would like to thank the referees for corrections and  helpful comments.

\end{document}